\documentclass[a4paper,12pt]{amsart}
\usepackage{epsfig,amsmath}

\begin{document}

\numberwithin{equation}{section}

\newtheorem{thm}[subsection]{Theorem}
\newtheorem{lem}[subsection]{Lemma}
\newtheorem{cor}[subsection]{Corollary}
\newtheorem{prop}[subsection]{Proposition}
\newtheorem{obs}[subsection]{Observation}

\theoremstyle{definition}
\newtheorem{definition}[subsection]{Definition}
\newtheorem{remark}[subsection]{Remark}

\newcommand{\R}{\mathbf{R}}
\newcommand{\US}{\mathbf{S}}

\def\a{\alpha}
\def\g{\gamma}
\def\s{\sigma}
\def\ts{\tilde\sigma}
\def\l{\lambda}
\def\t{\theta}
\def\Q{\mathcal{Q}}
\def\C{\mathcal{C}}
\def\cc{\mathrm{Cc}}
\def\eps{\varepsilon}
\def\scope{\mathrm{scope}}
\def\intr{\mathrm{int}}
\def\bd{\partial}
\def\0{\mathbf{0}}

\long\def\symbolfootnote[#1]#2{
	\begingroup
	\def\thefootnote{\fnsymbol{footnote}}\footnote[#1]{#2}
	\endgroup}

\parskip 5pt
\parindent 0pt
\baselineskip 15pt

\author{Bruce Solomon}
\address{Math Department, Indiana University, Bloomington IN 47405}
\email{solomon@indiana.edu}
\urladdr{mypage.iu.edu/$\sim$solomon}

%\subjclass{53A04}
\date{First draft December 2006. Last Typeset \today.}
\thanks{A fellowship from the Lady Davis foundation helped support this work}.

\begin{abstract}
We produce \emph{skew} loops---loops having no pair of parallel tangent lines---homotopic to any loop in a flat torus or other quotient of $\,\R^n\,$. The interesting case here is $\,n=3\,$.  More subtly for any $\,n\,$, we characterize the homotopy classes that will contain a skew loop having a specified loop $\,\tau\subset\US^{n-1}\,$ as tangent indicatrix.
\end{abstract}

\title{skew loops in flat tori.}
\maketitle

%%%%%%%%%%%%%%%%%% INTRODUCTION %%%%%%%%%%%%%%%%%%%

\section{Introduction} 
Call a differentiable curve in $\,\R^n$ \emph{skew} if it lacks any pair of parallel tangent lines. 

Skew curves and even skew loops are generic when $\,n>3\,$, but the basic examples in \cite{segre} and \cite{ghomi1} show that it takes ingenuity to find skew {loops} when $\,n=3\,$. They abound even then, however, existing in every knot-class \cite{wu}, or through any finite set of directed points \cite{ghomi2}. They are also interesting: Skew loops exist on all compact surfaces in $\,\R^3\,$ \emph{except} ellipsoids \cite{gs}, and thin tubes around skew loops have no disconnected shadows, disproving a conjecture that only ovaloids had this property \cite{ghomi1}, \cite{ghomi2}.

The definition of ``skew'' makes sense in $\,\R^{n}\,$ because translation identifies all tangent spaces in a canonical way. Tangent lines at distinct points of an arc are \emph{parallel} if they coincide under this identification.  But the identification persists when we descend from $\,\R^n\,$ to its quotient by a discrete subgroup $\,G\,$. So the definition of skew loop makes equally good sense in $\,\R^n/G\,$.

In particular, any skew loop in $\,\R^n\,$ descends to a homotopically trivial skew loop in $\,\R^n/G\,$. In \S\ref{sec:helices} below, we complement this fact by constructing explicit helical skew loops in each homotopically \emph{non-trivial} class of loops. Perhaps surprisingly, they are easier to construct than their homotopically trivial antecedents.

In \S\ref{sec:cc}, we determine, given a loop $\,\tau\,$ in the unit sphere $\,\US^{n-1}\,$, which homotopy classes in $\,\R^{n}/G\,$ admit skew loops having $\,\tau\,$ as \emph{tantrix}:

\begin{definition} The \emph{tantrix} of a $\,C^1\,$ curve $\,\s\,$ in $\R^n/G\,$ is the normalized velocity curve  $\,\tau\,$ mapping the domain of $\,\s\,$ to $\,\US^{n-1}\,$ via $\,\tau(t):=\dot\s/|\dot\s|\,$.\symbolfootnote[2]{As in \cite{s1}, \cite{s2} and \cite{gs}, \emph{tantrix} contracts the unwieldy traditional term \emph{tangent indicatrix} }\end{definition}

We leave the domain of $\,\s\,$  unspecified above because of a technical dilemma:  One defines \emph{$\,C^1\,$ loop} most elegantly as an immersion of the circle. But we also need to regard such a loop as an immersed \emph{arc} $\,\s:[a,b]\to\R^n/G\,$ with $\,\s(b)=\s(a)\,$ and $\,\dot\s(a)=\dot\s(b)\,$. We will move between these equivalent formulations without further comment, trusting the reader to understand that when, for instance, we call a loop \emph{embedded}, we have the first formulation in mind, while we prefer the second when we lift a loop in $\,\R^n/G\,$ to an arc in $\,\R^n\,$.

Either way, a loop $\,\s\,$ in $\,\R^n/G\,$ clearly has parallel tangent lines at inputs $\,t_1,t_2\in\US^1\,$ iff its tantrix $\,\tau\,$ satisfies $\,\tau(t_2) = \pm \tau(t_1)\,$. Hence
\begin{obs}\label{obs:crit}
A loop in $\,\R^n/G\,$ is skew iff its tantrix is embedded and disjoint from its own antipodal image.
\end{obs}

On the other hand, a loop $\,\tau\subset\US^{n-1}\,$ can satisfy both conditions above \emph{without} forming the tantrix of any loop in $\,\R^n\,$. This is the problem Segre, Ghomi and others had to solve in the papers mentioned above.  As the rank of $\,G\,$ increases, this problem gets easier, but a new problem arises: A loop in $\,\US^{n-1}\,$ may be embedded and avoid its antipodal image without forming the tantrix of any loop in a \emph{particular homotopy class of $\,\R^{n}/G\,$}.  After preparing some facts about convex cones, we settle this new issue in Theorem \ref{thm:main}.

\goodbreak
%%%%%%%%%%%%%%% HELICAL EXAMPLES %%%%%%%%%%%%%%%%%%%

\section{Examples: Helical skew loops of all homotopy types}\label{sec:helices}

Consider any discrete subgroup  $\,G\subset\R^n\,$, and write  $\,\pi:\R^n\to\R^n/G\,$ for the covering projection. Recall that one can lift any loop $\,\s:[a,b]\to\R^n/G\,$ to an arc $\,\a:[a,b]\to\R^n\,$ with $\,\pi\circ \a = \s\,$. Since $\,\s(b)=\s(a)\,$, we must then have $\,\a(b)-\a(a)= g \,$ for some $\, g \in G\,$.

Conversely, any {arc} $\,\a:[a,b]\to\R^{n}\,$ with $\,\a(b)-\a(a)= g \in G\,$ clearly projects to a loop in $\,\R^{n}/G\,$ homotopic to the geodesic  $\,t\mapsto \pi(t\, g )\,$, $\,a\le t\le b\,$. 

\begin{definition} 
Given an arc $\,\a:[a,b]\to\R^{n}\,$ with $\,\a(b)-\a(a)= g \in G\,$, we will call the loop $\,\pi\circ \a\,$ in $\,\R^{n}/G\,$ \emph{homotopic to $\, g \,$} or simply \emph{$g$-homotopic}, reflecting the familiar isomorphism $\,\pi_{1}\left(\R^{n}/G\right)\approx G\,$.
\end{definition}

\begin{prop}
Choose any non-zero $\, g \in G\,$, any orthonormal pair $\,u:=\{u_1,u_2\}\,$ perpendicular to $\, g \,$, and any $\, r>0\,$. Then the helical arc $\,h_{ g ,u}:[0,2\pi]\to\R^n\,$ given by
\[
h_{ g ,u}(t):= t\, g  +  r\left (u_1\,\cos t + u_2\,\sin  t\right)
\]
projects to a $g$-homotopic skew loop in $\,\R^{n}/G\,$.

\end{prop}

\begin{proof}
Clearly $\,h_{ g ,u}(2\pi)-h_{ g ,u}(0) =  g \,$, making $\,h_{ g ,u}\,$ homotopic to  $\, g \,$. To prove skewness, compute its tantrix:
\[
{h'_{ g ,u}(t)\over\left|h'_{ g ,u}(t)\right|}={ g  \over\sqrt{| g |^2 +  r^2}} + {r\over\sqrt{| g |^2 +  r^2}}\,\left (u_2\,\cos t - u_1\,\sin t\right)\ .
\]
This is clearly an embedded circle in $\,\US^{n-1}\,$ with radius $\, {r\big/\sqrt{| g |^2+ r^2}}\,$. Having radius less than 1, it avoids its own antipodal image, so by Observation \ref{obs:crit}, we are done.
\end{proof}

Given a non-zero $\, g \in G\,$, and letting $\,r\,$ vary between $\,0\,$ and $\,\infty\,$,  the tantrices above assume every radius between 0 and 1. Our helical examples thus make every small circle in  $\,\US^{n-1}\,$ centered at $\, g /| g |\,$ the tantrix of skew loop homotopic to $\, g \,$. We now want to generalize this fact, determining when an \emph{arbitrary} loop in $\,\US^{n-1}\,$ forms the tantrix of a $g$-homotopic loop in $\,\R^{n}/G\,$ for some pre-assigned $\,g\in G\,$.

%%%%%%%%%%%%%%%%% CONVEX CONES %%%%%%%%%%%%%%%%%%%%

\section{Convex cones and tantrices}\label{sec:cc}

A \emph{cone} in $\,\R^{n}\,$ is a subset invariant under all positive dilations. A set is \emph{convex} if it contains the line segment joining any pair of its points.

\begin{definition}[Convex cone]
Let $A\subset\R^{n}$. Denote the smallest con- vex cone containing $\,A\,$ by $\,\cc(A)\,$, and call it the \emph{convex cone over $\,A\,$}.
\end{definition}

\begin{remark}\label{rem:ch}
The convex cone over $\,A\,$ is simultaneously (1) The convex hull of the cone over $\,A\,$, and (2) The cone over the convex hull of $\,A\,$, the \emph{convex hull} being the intersection of all convex supersets.  One argues these facts easily using the separation property of convex sets. The following version of that property will help us prove our key Lemma \ref{prop:scope}.
\end{remark}

Denote the interior of a set $\,A\,$ by $\,\intr\,A\,$ and its boundary by $\,\bd A\,$.

\begin{lem}[Separation]\label{lem:sep}
Suppose $\,\C\subset\R^{n}\,$ is a convex cone, and $\,p\not\in\C\,$.  Then some hyperplane through $\,\mathbf{0}\,$ separates $\,p\,$ from $\,\intr\,C\,$. Specifically, for some unit vector $\,u\in\R^{n}\,$, we have
\[
u\cdot p\le 0< u\cdot x\quad\text{for all $\,x\in \intr\,\C\,$},
\]
with equality on the left iff $\,p\in\partial\C\,$.
\end{lem}

\begin{proof}
The interior of $\,\,\C\,$ is clearly convex, and it is well-known \cite[Chap. II]{valentine} that any point $\,p\,$ outside an open convex set can by separated from it by a hyperplane. More precisely, for some unit vector $\,u\in\R^{n}\,$, some scalar $\,a\,$, and for all $\,x\in \intr\,\C\,$, we have
\begin{equation}\label{eqn:a2}
u\cdot p\le a < u\cdot x\ .
\end{equation}
Since $\,\C\,$ is a cone, its closure contains $\,\mathbf{0}\,$. So $\,a\le 0\,$, and the inequalities above must hold with $\,p\,$ replaced by $\,\lambda p\,$ for any $\,\lambda>1\,$. Make that replacement, divide by $\,\lambda\,$, and recall that  $\,\C\,$ is a cone to deduce
\[
u\cdot p\le {a\over \lambda} < u\cdot x\quad\text{for all $\,x\in\intr\,\C\,$ and all $\,\lambda>1\,$.}
\]
Now let $\,\lambda\to+\infty\,$ to see that
\[
u\cdot p\le 0 \le u\cdot x\quad\text{for all $\,x\in\intr\,\C\,$.}
\]
But equality on the right for some $\,x\in\intr\,\C\,$, would imply $\,x-\eps u\in\intr\,\C\,$ for small enough $\,\eps>0\,$, contradicting the same inequality.
\end{proof}

\begin{definition}[Full arc]
Call an arc $\,\a:I\to\R^{n}\,$ \emph{full} if its image spans $\,\R^{n}\,$. Full arcs are of course generic, and by assuming $\,\a\,$ full, we simplify many statements below. Of course one can always remove the assumption by replacing  $\,\R^{n}\,$ with the span of $\,\a\,$.
\end{definition}

The following statement about full arcs is our main tool for determining which loops in $\,\US^{n-1}\,$ form tantrices of loops in $\,\R^{n}/G\,$. It constitutes an ``if and only if'' version of a result due to Ghomi \cite[Lemma 2.3]{ghomi2}.

\begin{lem}\label{prop:scope}
Suppose an arc $\a:I\to\R^{n}$ is full, and $\,p\in\R^{n}\,$. Then $\,p\in\intr\,\cc(\a)\,$  iff there exists a continuous $\,\mu:I\to(0,\infty)\,$ which is constant on $\,\partial I\,$, and satisfies
\begin{equation}\label{eq:mass}
p=\int_{I}\mu(t)\,\a(t)\ dt\ .
\end{equation} 
\end{lem}

\begin{proof} As cited above, Ghomi shows that such a $\,\mu\,$ exists when $\,p\in\intr\,\cc(\a)\,$. Actually, he argues this in $\,\R^{3}\,$, quoting a ``6-point'' result of Steinitz for convex sets in $\,\R^{3}\,$. But Steinitz's result includes the analogous ``$2n$-point'' statement for convex sets in $\,\R^{n}\,$ \cite[Thm. 3.13]{valentine}, and Ghomi's lemma generalizes  without difficulty if one uses that fact.

The converse is easier\footnote{The referee points out that \cite[Lemma 2.1]{ghomi3} yields an alternative proof of this converse.}; we prove it by contradiction, assuming some positive function $\,\mu\,$ satisfies (\ref{eq:mass}) with  $\,p\not\in\intr\,\cc(\a)\,$. By Lemma \ref{lem:sep} we would then have a unit vector $\,u\,$ with
\begin{equation*}\label{eq:sep2}
u\cdot p\le 0< u\cdot x\quad\text{for all $\,x\in\intr\,\cc(\a)\,$}.
\end{equation*}
Since $\,\a(I)\subset {\cc(\a)}\,$,  the right-hand inequality above gives $\, u\cdot\a(t) \ge 0\,$ for all $\,t\in I\,$. Playing that against the left-hand inequality and the fact that $\,\mu> 0\,$, however, we get
\begin{equation}\label{eq:sharp}
0\ge u\cdot p =\int_I \mu(t)\,\left(u\cdot \a(t)\right)\ dt\ge 0\ ,
\end{equation}
whence the integral must vanish, forcing $\,u\cdot\a\equiv 0\,$. This  puts $\,\a\,$ in a hyperplane, contradicting our fullness assumption.
\end{proof}

%\begin{remark}\label{rem:gromov}
%Prefiguring Ghomi's lemma \cite[Lemma 3.2]{ghomi2}, Gromov \cite[p. 168]{gromov} showed that a loop $\,\tau:\US^1\to\US^{n-1}\,$ forms the tantrix of a loop in $\,\R^n\,$ if its convex hull contains a ball centered at the origin. In view of Remark \ref{rem:ch}, this follows immediately from the Corollary above in the case where $\,G\,$ is trivial and $\,\g=\0\,$. 
%\end{remark}
%
%\begin{remark}
%The ``proper subspace'' hypothesis here doesn't seriously restrict the Lemma---it simply means we should take the smallest subspace containing $\,\a(I)\,$ to be $\,\R^n\,$. The convex cone over $\,\a\,$ then has non-empty interior and the Lemma applies.
%\end{remark}

\begin{cor}\label{cor:gamma}
Let $\,G\subset\R^{n}\,$ be a discrete subgroup. Then a full loop $\,\tau:\US^1\to\US^{n-1}\,$ forms the tantrix of a $g$-homotopic loop in $\,\R^{n}/G\,$ iff $\,  g \in\intr\,\cc(\tau)\,$.
\end{cor}

\begin{proof}
If $\,g\in\intr\,\cc(\tau)\,$, the Lemma yields a continuous $\,\mu:[0,1]\to(0,\infty)\,$ with $\,\mu(0)=\mu(1)\,$, whence the formula
\[
\a(s):=\int_0^s\mu(t)\,\tau(t)\ dt
\]
parametrizes a $\,C^{1}\,$ arc in $\,\R^{n}\,$ with tantrix $\,\tau\,$, with $\,\a(0)=0\,$, and with $\,\a(1)= g \,$. As discussed early in \S\ref{sec:helices},  such an arc projects to a $g$-homotopic loop in $\,\R^{n}/G\,$.

Conversely, any immersed $\,C^{1}\,$ loop in $\,\R^{n}/G\,$ which is homotopic to 
$\, g \,$  and has tantrix $\,\tau\,$, must lift to an arc $\,\a:[0,1]\to\R^{n}\,$ with $\,\a(1)-\a(0)= g \,$ and $\,\a'(0)=\a'(1)\,$. The fundamental theorem of calculus then gives
\[
 g =\a(1)-\a(0) = \int_0^{1} |\a'(t)|\,\tau(t)\ dt\ .
\]
Since $\,|\a'|\,$ is positive, continuous, and constant on $\,\partial[0,1]\,$, Lemma~\ref{prop:scope} ensures that $\, g \in\intr\,\cc(\a)\,$.
\end{proof}

When combined with Observation \ref{obs:crit}, Corollary \ref{cor:gamma} immediately tells us when a  loop in $\,\US^{n-1}\,$ forms the tantrix of a \emph{skew} $g$-homotopic loop: 

\begin{thm}\label{thm:main}
Suppose $\,G\subset\R^{n}\,$ is a discrete subgroup. Then a full loop $\,\tau:\US^{1}\to \US^{n-1}\,$ forms the tantrix of a $g$-homotopic skew loop in $\,\R^{n}/G\,$  iff $\,\tau\,$ 
\begin{enumerate}
\item 
is embedded, 
\item
avoids its antipodal image, and 
\item
contains $\, g\,$ in the interior of its convex cone.
\end{enumerate}
\end{thm}

\begin{cor} 
Suppose $\,G\subset\R^{n}\,$ is a discrete subgroup. If a full loop $\,\tau:\US^{1}\to\US^{n-1}\,$ forms the tantrix of a skew loop in $\,\R^{n}\,$, it forms the tantrix of a $g$-homotopic skew loop in $\,\R^{n}/G\,$ for \emph{every} $\,g\in G\,$.
\end{cor}

\begin{proof}
By Corollary \ref{cor:gamma}, $\,\tau\,$ is embedded, avoids its antipodal image, and satisfies $\,\0\in\intr\cc(\tau)\,$. But the latter implies $\,\cc(\tau)=\R^{n}\,$, so $\, \intr\cc(\tau)\,$ contains every $\,g\in G\,$.
\end{proof}

A full loop contained in a hemisphere can never be the tantrix of a loop (skew or not) in $\,\R^{n}\,$. This follows, e.g.,  from Corollary \ref{cor:gamma} with $\,G=\{\mathbf{0}\}\,$. But the situation in a flat torus is quite different:

\begin{cor}\label{cor:fullrank}
Suppose a discrete subgroup $\,G\subset\R^{n}\,$ has rank $\,n\,$. Then \emph{every} full embedded loop $\,\tau:\US^{1}\to\US^{n-1}\,$ that avoids its antipodal image is the tantrix of  a skew loop in the torus $\,\R^{n}/G\,$.
\end{cor}

\begin{proof}
Since $\,\tau\,$ is full, $\,\cc(\tau)\,$ contains the cone over a basis for $\,\R^{n}\,$. Sufficiently far from the origin, such a cone will contain arbitrarily large open balls. Since $\,G\,$ has rank $\,n\,$, all large enough balls contain elements of $\, G\,$. Theorem \ref{thm:main} then makes $\,\tau\,$ the tantrix of a skew loop in $\,\R^{n}/G\,$. \end{proof}

%\begin{definition}[Big and small arcs]
%We call a full arc $\,\tau\,$ in $\,\US^{n-1}\,$  \emph{small} if it lies in in a closed hemisphere; otherwise $\,\tau\,$ is \emph{big}. One easily sees that $\,\tau\,$ is big iff $\,\mathbf{0}\in\intr\,\cc(\tau)=\R^{n}\,$. 
%\end{definition}

We close with an observation about the most interesting case, $\,n=3\,$.  Suppose a loop $\,\tau\,$ bounds a region $\,\Omega\,$ compactly contained in an open hemisphere. Then $\,\tau\,$ automatically satisfies conditions (1) and (2) of Theorem \ref{thm:main}. Since any point in $\,\Omega\,$ lies in $\,\cc(\tau)\,$ (exercise) that Theorem now yields
\begin{cor}
Suppose $\,G\subset\R^{3}\,$ is a discrete subgroup, and $\,\tau:\US^{1}\to\US^{2}\,$ bounds a region $\,\Omega\,$ compactly contained in an open hemisphere. Then $\,\tau\,$ forms the tantrix of a $g$-homotopic skew loop in $\,\R^{3}/G\,$ whenever $\,g/|g|\in\Omega\,$.
\end{cor}

%%%%%%%%%%%%%%%  ACKNOWLEDGMENTS  %%%%%%%%%%%%%%%%%

\section*{Acknowledgments} We thank the the math departments at UCSD, and at the Technion in Haifa, Israel for their hospitality during this work. Fellowship support from the Lady Davis Foundation made our Technion stay possible.

%%%%%%%%%%%%%%%%%  BIBLIOGRAPHY %%%%%%%%%%%%%%%%%%%

\vfill

\end{document}